# A cluster identification framework illustrated by a filtering model for earthquake occurrences

ZHENGXIAO WU

*Department of Statistics and Applied Probability, National University of Singapore.*
*E-mail:* stawz@nus.edu.sg

A general dynamical cluster identification framework including both modeling and computation is developed. The earthquake declustering problem is studied to demonstrate how this framework applies.

A stochastic model is proposed for earthquake occurrences that considers the sequence of occurrences as composed of two parts: earthquake clusters and single earthquakes. We suggest that earthquake clusters contain a "mother quake" and her "offspring." Applying the filtering techniques, we use the solution of filtering equations as criteria for declustering. A procedure for calculating maximum likelihood estimations (MLE's) and the most likely cluster sequence is also presented.

## 1. Introduction

Suppose one observes a series of events $X_1, X_2, \ldots, X_n$ occurring at times $\tau_1, \tau_2, \ldots, \tau_n$. Each event is either "normal" or "abnormal." The objective is to identify those "abnormal" events.

One application of this problem is in epidemiology. For instance, the patients with Severe Acute Respiratory Syndrome (SARS) have symptoms similar to those of common flu patients. However since SARS is much more infectious than common flu, the SARS patients often appear in clusters. Such statistical evidence enables us to identify the SARS patients by mathematical tools. It provides a supplementary method to the costly medical test.

Another application is to collusion set detection. In a stock market, a group of traders forms a collusion set if they heavily trade among themselves in order to manipulate the stock price. It is of interest to catch this kind of malpractice as early as possible.

Considering each trade record as an event, it is intuitive that the malicious trading events tend to cluster. Assuming that a distance measuring the dissimilarity between any two records is available, Palshikar and Apte tackle the problem via graph clustering in [12]. They ignore the time stamp on the trade record so that a point process is reduced to a graph. But the temporal information is lost in their method.

These examples motivated our filtering model. We model the observations as a mixture of two independent marked point processes representing the "normal" and the "abnormal" events, respectively. Each new "abnormal" event will change the intensity of the "abnormal" point process. Typically the "abnormal" event increases the intensity for additional "abnormal" events in its neighborhood. Our goal is to compute the conditional probability of each observed event being abnormal in real time. Employing filtering techniques, we derive versions of the Zakai and Kushner–Stratonovich equations. Under a Markov condition, a sequential algorithm is presented to calculate the exact conditional probability that we are interested in.

Unfortunately, the data set for the two examples above is not available. We will present our methodology in the context of the "earthquake declustering problem." Even though there is no agreement on the underlying mechanism of earthquake occurrence in the current seismology literature, we want to emphasize that this example is mainly for the purpose of illustration. Our framework is for general modeling and computation. It could be adapted for different data sets in various areas.

It is well known that earthquakes often occur in clusters. The largest quake in a cluster is called the main shock, those before it are called foreshocks, and those after it are called aftershocks. The aftershocks in an earthquake swarm are relatively easy to predict. However, there are also many earthquakes that strike without any foreshocks or aftershocks. As the authors stated in [8]: "To forecast the location of the large earthquakes, it is necessary to analyze the background seismicity, for which removal of temporal cluster members is considered to be of central importance."

In this article, we propose a space–time point process model stemming from [14]. The observed earthquakes are considered as a mixture of earthquake swarms (a swarm contains at least two quakes) and single quakes. This could be considered as a special case of the "cluster processes" ([2], Section 6.3): a cluster process is composed of clusters that contain only a single point and clusters that have multiple points; in our model, we distinguish the single point events (single quakes) as the "noise" and the multiple point events (earthquake swarms) as the "signal." The conditional probability that a quake is in a cluster becomes a natural criterion for declustering. The filtering theory hence can be applied. We assume that, at most, one cluster is active at a time. This assumption can be relaxed with increased computational complexity.

In the literature, inference for partially observed stochastic processes is often obtained by using Markov chain Monte Carlo (MCMC) methods (see, e.g., [6]). A particle algorithm is also proposed in [14]. Such approximation methods are more flexible, but they are time-consuming and the approximation error is usually difficult to estimate. This paper deals with finding analytic solutions for some cases.

The paper is organized as follows: Section 2 describes the generic model and the filtering equations; Section 3 presents the computational procedure for the conditional expectation of interest under the "mother quake" assumption (the first quake in a cluster triggers all the other quakes in that cluster); Section 4 illustrates the numerical results for earthquakes in central and western Japan; Section 5 summarizes the conclusions and describes future work; Appendix A gives the algorithm that calculates the maximum likelihood estimators of the parameters and the most likely cluster path; finally, the proofs are contained in Appendix B.

## 2. The generic model

### 2.1. Formulation of the model

Suppose observed information about a quake is represented by a mark in a space $E$. For example, $E$ could be $\mathbf{R}^3$, recording the earthquake's magnitude and the epicenter's location longitude and latitude. We model observations as a marked point process $O$ with marks in $E$. $O$ is the mixture of two independent point processes $N$ and $C$, which stand for the single quakes and earthquake clusters, respectively. Hence letting $O(A,t)$ denote the number of quakes characterized by values in $A$ ($A$ is a subset of $E$) observed up to time $t$, we can write

$$O(A,t) = N(A,t) + C(A,t).$$

We assume that $N$ is a Poisson process with intensity $\gamma$ relative to a reference measure $\nu$, hence the single quake model is just a Poisson random measure on $E \times [0,\infty)$ with mean measure $\nu_0(\mathrm{d}u \times \mathrm{d}s) = \gamma(u,s)\nu(\mathrm{d}u)\,\mathrm{d}s$.

We model clusters to be randomly initiated and assume they eventually die out; we also assume that there is at most one active cluster at a time as mentioned in Section 1. Let $D$ be the process that indicates whether a cluster is active or not. The process $C$ adds a mark $u$ at time $s$ with non-negative predictable intensity $\lambda(u, s, D_{s-}, \eta_{s-})$, where $\eta$ is the configuration of both the marks and occurrence times of all the previous cluster quakes. More precisely, if cluster quake $c_i$ occurs at $t_i$, then $\eta_t = \sum_{\{i:\ t_i \leq t\}} \delta_{(c_i, t_i)}$, where $\delta_{(c_i, t_i)}$ is the Dirac measure concentrated on the point $(c_i, t_i)$. Therefore, $\eta$ is a counting measure on $E \times [0, \infty)$.

When $D=0$, there is no active cluster and an intensity $\lambda(u, s, 0, \eta_{s-})$ gives the rate at which a new cluster is initiated by an event with mark $u$ at time $s$. Once initiated, the cluster grows with intensity $\lambda(u, s, 1, \eta_{s-})$ until it dies out.

Under very mild conditions on the intensities (see [5]), the point processes can be written as solutions of stochastic differential equations. In particular, we can write

$$\begin{aligned}
O(A,t) &= N(A,t) + C(A,t) \\
&= \int_{A \times [0,\infty) \times [0,t]} \mathbf{1}_{[0,\gamma(u,s)]}(v) \xi_1(\mathrm{d}u \times \mathrm{d}v \times \mathrm{d}s) \\
&\quad + \int_{A \times [0,\infty) \times [0,t]} \mathbf{1}_{[0,\lambda(u,s,D_{s-},\eta_{s-})]}(v) \xi_2(\mathrm{d}u \times \mathrm{d}v \times \mathrm{d}s),
\end{aligned}$$

where $\xi_1$ and $\xi_2$ are independent copies of a Poisson random measure on $E \times [0, \infty) \times [0, \infty)$ with mean measure $\nu \times \ell \times \ell$, denoting Lebesgue measure by $\ell$.

In this article, we define $D$ as follows: $D$ is equal to 1 once a cluster is initiated; $D$ has a probability $p$ to die out (i.e., $D = 0$) whenever a new observation is added to the cluster; $D$ is independent of all previous history. Thus, for an arbitrary function $f(D_t, \eta_t)$,

$$f(D_t, \eta_t) = \int_{E \times [0,t]} \mathbf{1}_{\{D_{s-}=0\}}[f(1, \eta_{s-} + \delta_{(u,s)}) - f(0, \eta_{s-})] \\
+ \mathbf{1}_{\{D_{s-}=1\}}[f(1 - I_{C(E,s)}, \eta_{s-} + \delta_{(u,s)}) - f(1, \eta_{s-})]C(\mathrm{d}u \times \mathrm{d}s), \quad (2.1)$$

where the $\{I_k, k = 1, 2, \ldots\}$ are independent Bernoulli random variables with parameter $p$ that are also independent of $N$ and $C$. This follows by writing the right-hand side as a finite sum where most terms cancel out.

In practice, $f(D_t, \eta_t)$ contains information about $D_t$ and $\eta_t$. Statistical inferences can be drawn if we are able to compute the conditional expectation of $f$ based on the observations $O$. The rest of the paper mostly deals with how to realize such a computation for arbitrary $f$.

It is worth noting that our whole problem is essentially discrete and finite, hence the measurability of functions is (and should be) of minor concern. As $D_t$ is either 0 or 1, and $\eta_t$ can only take finitely many values as well ($2^n$ if there are $n$ observations), thus the function domain of $f$ is finite. Therefore, all the functions are measurable.

### 2.2. The filtering equations

We derive the filtering equation for the conditional distribution of $\eta$ given observations of $O$ using a reference measure approach. If $(\Omega, \mathcal{F}, Q)$ is a probability space and $P$ is a second probability measure on $\mathcal{F}$ given by $\mathrm{d}P = L\,\mathrm{d}Q$, then for any sub-$\sigma$-algebra $\mathcal{D} \subset \mathcal{F}$ and $L^1$-random variable $Z$,

$$E^P[Z|\mathcal{D}] = \frac{E^Q[ZL|\mathcal{D}]}{E^Q[L|\mathcal{D}]}.$$

We are going to use a reference probability measure $Q$ under which the observations have a relatively simple structure. In the following lemma, $N$ and $C$ are independent Poisson random measures under $Q$. We first introduce a definition that is used in the lemma.

**Definition 2.1.** *A Poisson process $N$ is **compatible** with a filtration $\{\mathcal{F}_t\}$ if $N$ is $\{\mathcal{F}_t\}$-adapted and $N(t + \cdot) - N(t)$ is independent of $\{\mathcal{F}_t\}$ for every $t \geq 0$.*

**Lemma 2.2.** *On $(\Omega, \mathcal{F}, Q)$, let $N$ and $C$ be independent Poisson random measures with mean measures $\nu_0(\mathrm{d}u \times \mathrm{d}s) = \gamma(u,s)\nu(\mathrm{d}u)\,\mathrm{d}s$ and $\nu_1(\mathrm{d}u \times \mathrm{d}s) = \lambda_Q(u,s)\nu(\mathrm{d}u)\,\mathrm{d}s$, respectively; let $D$ be a cadlag process independent of $N$. Assume all processes are compatible* with $\{\mathcal{F}_t\}$. $L$ is determined by solving

$$L(t) = 1 + \int_{E \times [0,t]} \left( \frac{\lambda(u, s, D_{s-}, \eta_{s-})}{\lambda_Q(u, s)} - 1 \right) L(s-)[C(\mathrm{d}u \times \mathrm{d}s) - \lambda_Q(u,s)\nu(\mathrm{d}u)\,\mathrm{d}s] \quad (2.2)$$

*and assuming that $L$ is a $\{\mathcal{F}_t\}$-martingale. Let $P$ satisfy $\mathrm{d}P_{|\mathcal{F}_t} = L(t)\,\mathrm{d}Q_{|\mathcal{F}_t}$. Then $P$ is a probability measure and under $P$, for all $A$ such that $\int_0^t \int_A \lambda(u,s,D_s,\eta_s)\nu(\mathrm{d}u)\,\mathrm{d}s < \infty$ for each $t > 0$,*

$$C(A,t) - \int_{A \times [0,t]} \lambda(u, s, D_s, \eta_s)\nu(\mathrm{d}u)\,\mathrm{d}s$$

*is a local martingale and $N$ is independent of $C$ and is a Poisson random measure with mean measure $\nu_0$.*

Thus under $P$ both $N$ and $C$ have the intensity described in Section 2.1. Hence Lemma 2.2 gives the form of the Radon–Nikodym derivative (or the likelihood) of $P$ with respect to $Q$. Our further computation then can be justified by the uniqueness of the martingale problem (see [9] or [4], Chapter 4). The assumption that $L$ is a $\{\mathcal{F}_t\}$-martingale is very mild.

*Remark 2.3.* The following condition is sufficient to ensure that (2.2) is a well-posed equation and that $L$ is a martingale.

(Condition 1) $\quad \nu(E) < \infty, \lambda_Q(u,s) \quad$ and $\quad \dfrac{\lambda(u,s,D_{s-},\eta_{s-})}{\lambda_Q(u,s)}$ are uniformly bounded.

The process $D$ has finitely many jumps in bounded time intervals. Thus we can record the history of the process $D$ by a counting measure $h_t = \sum_{\{i:\, t_i \leq t\}} \delta_{(D_{t_i}, t_i)}$; the sum is over those $t_i$ when $D$ takes jumps. Hence it represents a path that has value $D_{t_i}$ in time interval $[t_i, t_{i+1})$. As in (2.1), let $f$ be an arbitrary function on the two counting measures $(h_s, \eta_s)$, and set

$$\phi(f, s) = E^Q[f(h_s, \eta_s) L(s) | \mathcal{F}_s^O], \tag{2.3}$$

$$\pi(f, s) = E^P[f(h_s, \eta_s) | \mathcal{F}_s^O] = \frac{E^Q[f(h_s, \eta_s) L(s) | \mathcal{F}_s^O]}{E^Q[L(s) | \mathcal{F}_s^O]} = \frac{\phi(f, s)}{\phi(1, s)}. \tag{2.4}$$

Since $h_t$ contains all the information on $D_t$, we can write $D_t = D_t(h_t)$. Further, we abuse the notation a little and write $\lambda(u, s, D_{s-}(h_{s-}), \eta_{s-}) = \lambda(u, s, h_{s-}, \eta_{s-})$. We need this expression to simplify the notation in the following theorem and in the application in Section 4.

Let $\alpha$ denote the indicator that specifies whether or not a cluster is currently active, that is, $\alpha(h_s, \eta_s) = \mathbf{1}_{\{D_s=1\}} = D_s$. Let $q = 1 - p$ and

$$\begin{aligned}f_{\text{new}} &= [1 - \alpha(\cdot, \cdot)] f(\cdot + \delta_{(1,s)}, \cdot + \delta_{(u,s)}) \\ &\quad + \alpha(\cdot, \cdot)[p f(\cdot + \delta_{(0,s)}, \cdot + \delta_{(u,s)}) + q f(\cdot + \delta_{(1,s)}, \cdot + \delta_{(u,s)})].\end{aligned} \tag{2.5}$$

**Theorem 2.4.** *For an arbitrary function $f$ on $(h_s, \eta_s)$, let $\phi$, $\pi$ and $f_{\text{new}}$ be defined as in equations* (2.3)–(2.5). *Then $\phi$ and $\pi$ satisfy the stochastic integral equations*

$$\phi(f,t) = \phi(f,0) - \int_{E\times[0,t]} \phi(f(\cdot,\cdot)[\lambda(u,s,\cdot,\cdot) - \lambda_Q(u,s)], s)\nu(\mathrm{d}u)\,\mathrm{d}s$$

$$+ \int_{E\times[0,t]} \phi\left(f_{\text{new}}\frac{\lambda(u,s,\cdot,\cdot)}{\lambda_Q(u,s)} - f(\cdot,\cdot), s-\right)\frac{\lambda_Q(u,s)}{\lambda_Q(u,s) + \gamma(u,s)}O(\mathrm{d}u \times \mathrm{d}s)$$

*and*

$$\pi(f,t) = \pi(f,0)$$

$$+ \int_{E\times[0,t]} \frac{\pi(f_{\text{new}}\lambda(u,s,\cdot,\cdot), s-) - \pi(\lambda(u,s,\cdot,\cdot), s-)\pi(f,s-)}{\pi(\lambda(u,s,\cdot,\cdot) + \gamma(u,s), s-)} O(\mathrm{d}u \times \mathrm{d}s)$$

$$- \int_{E\times[0,t]} \{\pi(f(\cdot,\cdot)\lambda(u,s,\cdot,\cdot), s) - \pi(f,s)\pi(\lambda(u,s,\cdot,\cdot), s)\}\nu(\mathrm{d}u)\,\mathrm{d}s.$$

In Section 4, we will take $f$ as the indicator functions that indicate the status of $\eta$, so that $\pi(f,t)$ gives us the conditional probability that an observation is in the cluster.

## 3. Solutions of the filtering equation

Unlike the infinite-dimensional nonlinear filtering problem, the solution of which can only be approximated, the function space in our problem allows a natural finite decomposition since we have a finite function domain, that is, all the possible combinations of each observed event being in a quake swarm or not. Thus the exact solution could be computed theoretically, but generally the computational load increases exponentially as the number of observations increases. That is not feasible for online updating.

In this section and also in Appendix A, we assume that when a cluster is active, the cluster is assumed to be triggered by the first quake (mother quake); when no cluster is active, a new cluster will be initiated randomly with an intensity $\varepsilon$. To be precise, suppose one observes $u_i$ at time $\tau_i$. Let $y_i = (u_i, \tau_i)$ and the set of observations by time t be $O(t) = \{y_1, y_2, \ldots, y_k : \tau_k \leq t < \tau_{k+1}\}$. Then we have

$$\lambda(u, t, D_{t-}, \eta_{t-}) = D_{t-} \sum_{i=1}^{k} \lambda(u, t, y_i)\theta_0(y_i) + (1 - D_{t-})\varepsilon(u, t), \tag{3.1}$$

where $\theta_0(y_i) = \mathbf{1}_{\{y_i \text{ is the mother quake in the currently active cluster}\}}$ and $\theta_0(y_i)$ is defined as 0 if there is no active cluster at that time. We suppose that the functional form of $\lambda(u, t, y_i)$ is known. For the application in Section 4, $\lambda(u, t, y_i)$ is a Gaussian kernel (4.1) that does not depend on $t$. Note that there is, at most, one $\theta_0(y_i)$ ($i = 1, 2, \ldots, k$) non-zero at any moment $t$. The simple fact that $\theta_0(y_i)\theta_0(y_j) = 0$ if $i \neq j$ makes finding an analytic solution possible (see the proof of the following theorems). Formally, we can think of the intensity $\lambda$ as a vector with component $\lambda(u,t,y_i)$ at each "orthogonal" direction $\theta_0(y_i)$, $i = 1, 2, \ldots, k$. With the help of this kind of "orthogonal decomposition" of the function space, the problem can be reduced to be of polynomial complexity.

For simplicity, we also assume that there is no cluster active at time 0. Define $a(y,t) = \int_E \lambda(u,t,y)\nu(\mathrm{d}u)$ for $y \in E$ and $\varepsilon(t) = \int_E \varepsilon(u,t)\nu(\mathrm{d}u)$. The following two theorems give the algorithm to compute $\pi(f,t)$, $f$ is an arbitrary function of $D_s$ and $\eta_s$. Recall that $\alpha(D_s, \eta_s) = \mathbf{1}_{\{D_s=1\}} = D_s$.

**Theorem 3.1.** *For $\tau_k \leq t < \tau_{k+1}$*

$$\pi(\theta_0(y_i)\alpha, t) = \pi(\theta_0(y_i)\alpha, \tau_k) \mathrm{e}^{-\int_{\tau_k}^{t} a(y_i,s) - \varepsilon(s)\,\mathrm{d}s} b_k(t), \tag{3.2}$$

*where*

$$b_k(t) = \frac{1}{\sum_{j=1}^{k} \pi(\theta_0(y_j)\alpha, \tau_k) \mathrm{e}^{-\int_{\tau_k}^{t} a(y_j,s) - \varepsilon(s)\,\mathrm{d}s} + 1 - \sum_{j=1}^{k} \pi(\theta_0(y_j)\alpha, \tau_k)}$$

*and*

$$\pi(\theta_0(y_i)\alpha, \tau_{k+1}) = \frac{\pi(\theta_0(y)\alpha, \tau_{k+1}-)(\gamma_{k+1} + q\lambda_{k+1,i})}{d_{k+1}}, \qquad i < k+1,$$

$$\pi(\theta_0(y_{k+1})\alpha, \tau_{k+1}) = \frac{\sum_{j=1}^{k} (q\lambda_{k+1,j} - \varepsilon_{k+1}) \pi(\theta_0(y_j)\alpha, \tau_{k+1}-) + \varepsilon_{k+1}}{d_{k+1}},$$

*where $\gamma_{k+1} = \gamma(y_{k+1}, \tau_{k+1})$, $\lambda_{k+1,i} = \lambda(y_{k+1}, \tau_{k+1}, y_i)$, $\varepsilon_{k+1} = \varepsilon(y_{k+1}, \tau_{k+1})$ and $d_{k+1} = \sum_{j=1}^{k} (\lambda_{k+1,j} - \varepsilon_{k+1}) \pi(\theta_0(y_j)\alpha, \tau_{k+1}-) + \varepsilon_{k+1} + \gamma_{k+1}$.*

In Theorem 3.3, we can solve for $\pi(\theta_0(y_i)\alpha, t)$ as the first step in our algorithm. The task of computing $\pi(f,t)$ for more general $f$ is completed in the next theorem. Note that the solution $\pi(\theta_0(y_i)\alpha, t)$ is needed in (3.4).

**Theorem 3.2.** *For $\tau_k \leq t < \tau_{k+1}$,*

$$\pi(\theta_0(y_i)\alpha f, t) = \pi(\theta_0(y_i)\alpha f, \tau_k) \mathrm{e}^{-\int_{\tau_k}^{t} a(y_i,s) - \varepsilon(s)\,\mathrm{d}s} b_k(t), \tag{3.3}$$

$\pi(f,t)$ *satisfies*

$$\frac{\mathrm{d}\pi(f,t)}{\mathrm{d}t} = -\sum_{x \in Y(t)} \pi(\theta_0(x)\alpha f, t)[a(x,t) - \varepsilon(t)] + \pi(f,t) \sum_{x \in Y(\tau_k)} \pi(\theta_0(x)\alpha, t)[a(x,t) - \varepsilon(t)], \tag{3.4}$$

and

$$\pi(\theta_0(y_i)\alpha f, \tau_{k+1})$$
$$= \frac{\pi(\theta_0(y)\alpha f, \tau_{k+1}-)\gamma_{k+1} + q\lambda_{k+1,i}\pi(f(1,\cdot+\delta_{y_{k+1}})\theta_0(y_j)\alpha, \tau_{k+1}-)}{d_{k+1}}, \qquad i < k+1,$$
$$\pi(\theta_0(y_{k+1})\alpha f, \tau_{k+1}) = \frac{\pi(f(1,\cdot+\delta_{y_{k+1}})(1-\alpha), \tau_{k+1}-)\varepsilon_{k+1}}{d_{k+1}},$$
$$\pi(f, \tau_{k+1}) = \frac{\pi(f, \tau_{k+1}-)\gamma_{k+1} + \pi(f(1,\cdot+\delta_{y_{k+1}})(1-\alpha)\varepsilon_{k+1}, \tau_{k+1}-)}{d_{k+1}}$$
$$+ \frac{\sum_{j=1}^{k}\lambda_{k+1,j}\pi([pf(0,\cdot+\delta_{y_{k+1}}) + qf(1,\cdot+\delta_{y_{k+1}})]\theta_0(y_j)\alpha, \tau_{k+1}-)}{d_{k+1}},$$

where $f(0,\cdot+\delta_{y_{k+1}}) = f(\cdot+\delta_{(0,\tau_{k+1})},\cdot+\delta_{y_{k+1}})$, $f(1,\cdot+\delta_{y_{k+1}}) = f(\cdot+\delta_{(1,\tau_{k+1})},\cdot+\delta_{y_{k+1}})$.

Combining Theorems 3.3 and 3.2, we can compute $\pi(f,t)$ for an arbitrary $f$ in real time.

## 4. Application to an earthquake data set

We use the same data set as in [8]: the earthquakes in the period of 1926–1995 in the rectangular area 34°–39°N and 131°–140°E with magnitudes greater than 4.0 and depths less than 100 km.

We take $\nu$ to be the uniform measure on the rectangular region, $\gamma(u) = \gamma$ and $\lambda(u, D_t, \eta_t) = \mathbf{1}_{\{D_t=1\}}\sum_{i=1}^{k}\lambda(u,y_i)\theta_0(y_i) + \mathbf{1}_{\{D_t=0\}}\varepsilon$, where $\lambda(u,y_i)$ is proportional to a bivariate normal density:

$$\lambda(u, y_i) = \lambda \exp(-\|u - u_i\|^2/2d)/(2\pi d). \tag{4.1}$$

For our mother quake model, the maximum likelihood estimations (MLE's) are $\widehat{\gamma} = 0.1070$, $\widehat{\lambda} = 1.3274$, $\widehat{\varepsilon} = 0.0126$, $\widehat{d} = 0.0070$, $\widehat{p} = 0.2035$. The log-likelihood is $\log(L) = -21604$.

We are interested in

$$\theta(y)(\cdot,\cdot) = \mathbf{1}_{\{y \text{ is a quake in a cluster}\}}$$

and

$$D_t(\cdot,\cdot) = \mathbf{1}_{\{a \text{ cluster is active at time } t\}}.$$

We compute $\pi(\theta(y_i)(\cdot,\cdot),T)$ for all observed $y_i$ according to Theorem 3.2, where $T$ is the last moment of the year 1995. The results are compared with [8], where the authors declustered the observations by computing aftershock probabilities under an ETAS model. In Figure 1, plot (a) is the histogram of the aftershock probabilities as presented in [8]. We denote their aftershock probabilities as $p_1$. Plot (b) shows the distribution of the conditional probabilities $p_2$ in the mother quake model. We are pleased to see that our stochastic models give relatively deterministic answers. Around 95% of quakes have a probability of being in clusters that is either smaller than 0.1 or greater than 0.9, as can be seen in plot (b).

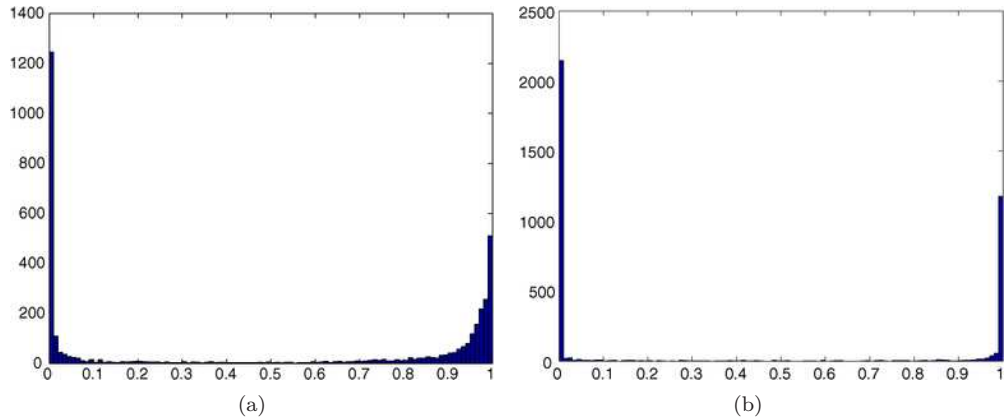

**Figure 1.** Histograms of (a) aftershock probabilities under ETAS model and (b) conditional probabilities to be in the cluster process in mother quake model.

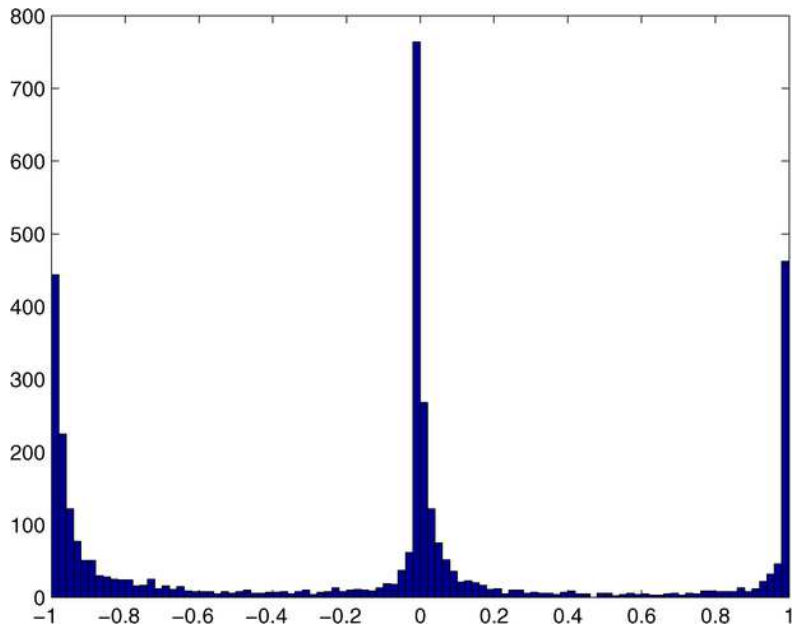

**Figure 2.** Histogram of $p_2 - p_1$.

Although both results have a bimodal shape, the one in [8] disagrees with our models for many individual quakes. This can be seen from Figure 2. The histogram presents the difference of these two probabilities.

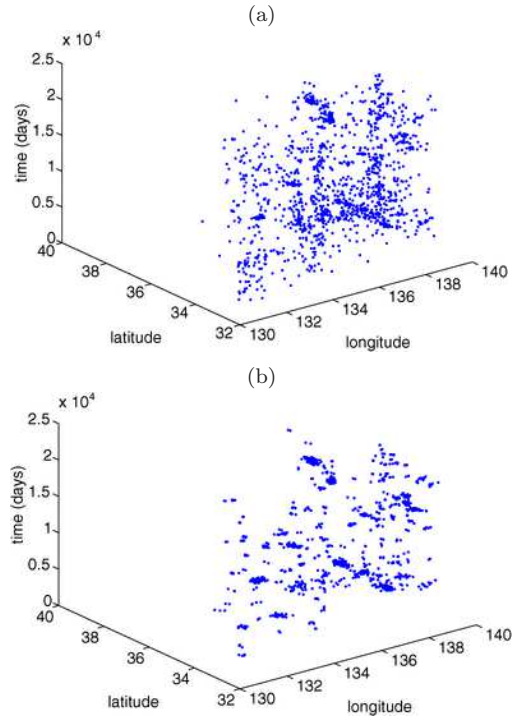

**Figure 3.** Time–space plots of the 1500 likely clustered earthquakes under different models: (a) under ETAS model; (b) under the mother quake model.

It seems that the data set supports our model more. We plot the earthquake clusters in each setting by removing quakes with a low probability of being in a cluster. The time-space plots in Figure 3 have 1500 quakes. The vertical axis represents time (unit in days). It is quite clear that the plots from our models have a stronger cluster pattern. The three-dimensional plots are available from http://www.stat.nus.edu.sg/~stawz/ and can be rotated and viewed in different perspectives.

We also can compute $\pi(D_{\tau_i}(\cdot,\cdot),T)$ to see the status of the cluster at different times. Under the mother quake assumption, Figure 4(a) gives us the conditional probability that the earthquake cluster is active. The answer is again quite distinct. Figure 4(b) shows that most conditional probabilities are either close to 0 or 1.

## 5. Discussion

Assumption (3.1) is just an example. Another earthquake model called the "domino" model is given in [14], where we assume that the last quake triggers the next quake in the cluster. It turns out that the mother quake model is more likely in our data set by comparing their likelihood. Roughly speaking, as long as the conditional intensity $\lambda$ modeling the cluster only depends on a "small" portion of the history (in (3.1), it only depends on the last mother quake), we can adopt an "orthogonal decomposition" and find an algorithm to find the analytic solution.

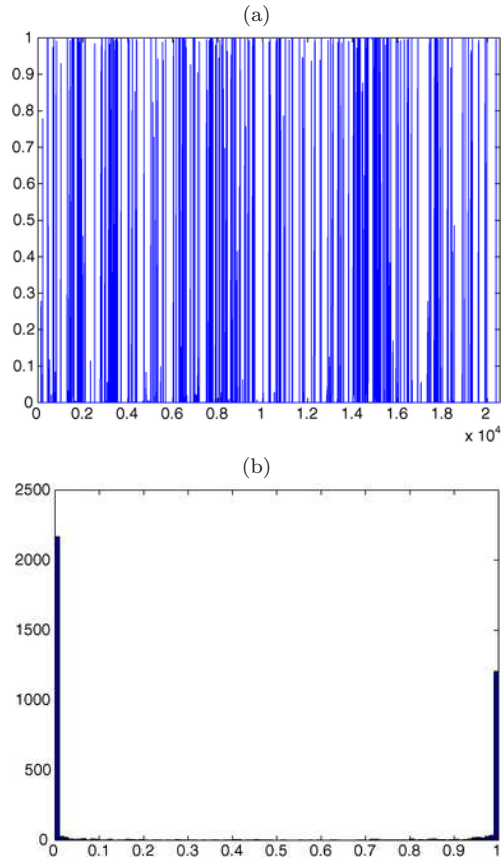

**Figure 4.** Under the mother quake assumption: (a) the conditional probability that the cluster was alive vs. time; (b) histogram of the conditional probability.

Thus assuming that, at most, one cluster is active at a time is not essential for our method. This simplified assumption prevents the presentation from getting more messy. We can similarly work out a decomposition if we assume that, at most, say, three clusters are active at a time.

Our filter separates the data set into the cluster quakes and the single quakes. Further data analysis in [14] shows geophysical differences. In particular, the magnitude of the cluster quakes is significantly different from the single quakes. The mother quakes are also significantly bigger than the offspring quakes. Note that we did not incorporate the magnitudes of the earthquakes into the model. This surprising finding further supports our model.

The application to seismology is only a special case of the filtering approach to abnormal cluster identification proposed in [14]. Other possible applications include epidemiology, intrusion detection in network security, criminology and quality control.

# Appendix A: Likelihood and maximum likelihood estimators

## A.1. Likelihood

**Theorem A.1.** *Let $\nu$ be a finite measure and, on $(\Omega, \mathcal{F}, Q)$, let $N$ and $C$ be independent Poisson random measures with mean measures $\nu(\mathrm{d}u)\,\mathrm{d}s$; let $D$ satisfy (2.1). Assume all processes are compatible with $\{\mathcal{F}_t\}$. Define $L_N$ and $L_C$ by solving*

$$L_N(t) = 1 + \int_{E \times [0,t]} [\gamma(u,s) - 1] L_N(s-)[N(\mathrm{d}u \times \mathrm{d}s) - \nu(\mathrm{d}u)\,\mathrm{d}s], \qquad (A.1)$$

$$L_C(t) = 1 + \int_{E \times [0,t]} [\lambda(u,s,D_{s-},\eta_{s-}) - 1] L_C(s-)[C(\mathrm{d}u \times \mathrm{d}s) - \nu(\mathrm{d}u)\,\mathrm{d}s] \quad (A.2)$$

*and assume that they are $\{\mathcal{F}_t\}$-martingales. Let $L = L_N L_C$. $L$ will also be an $\{\mathcal{F}_t\}$-martingale. Let $P$ satisfy $\mathrm{d}P_{|\mathcal{F}_t} = L(t)\,\mathrm{d}Q_{|\mathcal{F}_t}$. Then $P$ is a probability measure and under $P$, for all $A$ such that $\int_0^t \int_A \lambda(u,s,D_s,\eta_s)\nu(\mathrm{d}u)\,\mathrm{d}s < \infty$ for each $t > 0$,*

$$C(A,t) - \int_{A \times [0,t]} \lambda(u,s,D_s,\eta_s)\nu(\mathrm{d}u)\,\mathrm{d}s$$

*is a local martingale and $N$ is independent of $C$ and is a Poisson random measure with mean measure $\gamma(u,s)\nu(\mathrm{d}u)\,\mathrm{d}s$.*

*Remark A.2.* The $L$ derived from the theorem is the likelihood of our observation, which is the mixture of two processes. It is necessary to have it for the estimation of parameters. While in Lemma 2.2, the simplified version (2.2) is sufficient for proving Theorem 2.4, since it only concerns $f(D_t, \eta_t)$, which does not involve the process $N$. By applying $L$ derived here, we can prove a more general form of Theorem 2.4 so that we can have filtering equations about $f(N_t, D_t, \eta_t)$. We omit it because the notation gets worse and we do not use it in our application.

Our goal is to compute $E^Q[L|\mathcal{F}^O]$, the likelihood in our model. We can first solve (A.1) and (A.2):

$$L_N(t) = \exp\biggl\{\int_{E \times [0,t]} \log \gamma(u,s) N(\mathrm{d}u \times \mathrm{d}s)$$
$$- \int_{E \times [0,t]} [\gamma(u,s) - 1]\nu(\mathrm{d}u)\,\mathrm{d}s\biggr\},$$

$$L_C(t) = \exp\left\{\int_{E\times[0,t]} \log \lambda(u,s,D_{s-},\eta_{s-})C(\mathrm{d}u\times\mathrm{d}s)\right.$$
$$\left. - \int_{E\times[0,t]} [\lambda(u,s,D_{s-},\eta_{s-}) - 1]\nu(\mathrm{d}u)\,\mathrm{d}s\right\}.$$

Let

$$C(A,t) = \int_{A\times[0,t]} \rho_{O(E,s)}(u) O(\mathrm{d}u\times\mathrm{d}s), \tag{A.3}$$

where $\rho$ is the indicator of whether the observation is a cluster point. Under reference measure $Q$, $\rho_1,\rho_2,\ldots$ are i.i.d. Bernoulli$(1/2)$, and

$$L(t) = L_N(t)L_C(t)$$
$$= \exp\left\{\int_{E\times[0,t]} [(1-\rho_{O(E,s)})\log\gamma(u,s) + \rho_{O(E,s)}\log\lambda(u,s,D_{s-},\eta_{s-})]O(\mathrm{d}u\times\mathrm{d}s)\right.$$
$$\left. - \int_{E\times[0,t]} [\lambda(u,s,D_{s-},\eta_{s-}) + \gamma(u,s) - 2]\nu(\mathrm{d}u)\,\mathrm{d}s\right\}$$
$$\propto \exp\left\{\int_{E\times[0,t]} [(1-\rho_{O(E,s)})\log\gamma(u,s) + \rho_{O(E,s)}\log\lambda(u,s,D_{s-},\eta_{s-})]O(\mathrm{d}u\times\mathrm{d}s)\right.$$
$$\left. - \int_{E\times[0,t]} (\lambda(u,s,D_{s-},\eta_{s-}) + \gamma(u,s))\nu(\mathrm{d}u)\,\mathrm{d}s\right\}.$$

We use "$\propto$" since we ignore a constant, which has no impact on likelihood inference.

Recall $\gamma_j = \gamma(y_j,\tau_j)$, $\lambda_{j,i} = \lambda(y_j,\tau_j,y_i)$, $\varepsilon_j = \varepsilon(y_j,\tau_j)$. Let $\tau_0 = 0$ and suppose one observes $u_i$ at time $\tau_i$, $y_i = (u_i,\tau_i)$.

$$L(\tau_n) \propto \exp\left\{-\int_{E\times[0,\tau_n]} \gamma(u,s)\nu(\mathrm{d}u)\,\mathrm{d}s\right\}$$
$$\times \prod_{i=1}^n \gamma_i^{1-\rho_i} \lambda(u_i,\tau_i,D_{\tau_i-},\eta_{\tau_i-})^{\rho_i} \exp\left\{-\int_{E\times[\tau_{i-1},\tau_i]} \lambda(u,s,D_{s-},\eta_{s-})\nu(\mathrm{d}u)\,\mathrm{d}s\right\}. \tag{A.4}$$

Note that, conditional on observations $O$, all the possible partitions of $O = C + N$ are equally likely under $Q$. The conditional distribution of $L_N L_C$ is completely determined by $\{\rho_i\}$ and $\{I_k\}$. $\{I_k\}$ is independent of $\{\rho_i\}$ since it is independent of $N$ and $C$. Hence the computation of $E^Q[L|\mathcal{F}^O]$ is reduced to calculating the expectation of $L(\{\rho_i\},\{I_k\})$ with i.i.d. Bernoulli$(1/2)$ $\{\rho_i\}$ and i.i.d. Bernoulli$(p)$ $\{I_k\}$. This expectation can be expressed as a weighted sum over all possible values of $\{\rho_i\},\{I_k\}$.

We recall the special form of $\lambda(u,t,D_t,\eta_t)$:

$$\lambda(u_i,\tau_i,D_{\tau_i-},\eta_{\tau_i-}) = D_{\tau_i-}\sum_{j=1}^{i-1}\lambda_{i,j}\theta_{0i}(y_j) + (1-D_{\tau_i-})\varepsilon_i,$$

where $\theta_{0i}(y_j) = \mathbf{1}_{\{y_j \text{ is the latest mother quake in the cluster at time } \tau_i-\}}$ and $\lambda_{i,j}$ is defined as in Theorem 3.3. Let $E_i$ be the index of the latest mother quake at time $\tau_i-$, then $\lambda(u_i, \tau_i, D_{\tau_i-}, \eta_{\tau_i-}) = D_{\tau_i-}\lambda_{i,E_i} + (1 - D_{\tau_i-})\varepsilon_i$.

Assume no cluster is active at time 0, hence $D_0 = 0$. Therefore,

$$L(\tau_1) \propto \gamma_1^{1-\rho_1} \varepsilon_1^{\rho_1} \exp\left\{-\int_{E \times [0,\tau_1]} [\varepsilon(u,s) + \gamma(u,s)]\nu(\mathrm{d}u)\,\mathrm{d}s\right\}.$$

Sum over two terms corresponding to $\rho_1 = 0$ and $\rho_1 = 1$, respectively, and we have

$$E^Q[L(\tau_1)|\mathcal{F}^O] \propto \frac{(\varepsilon_1 + \gamma_1)}{2} \exp\left\{-\int_{E \times [0,\tau_1]} [\varepsilon(u,s) + \gamma(u,s)]\nu(\mathrm{d}u)\,\mathrm{d}s\right\}.$$

### A.2. The forward algorithm and MLE

It is not practical to sum over all the terms by brute force since the number of terms increases exponentially with respect to the number of observations. However, we can reduce the complexity to $O(n^2)$ by computing $\exp\{\int_{E \times [0,\tau_n]} \gamma(u,s)\nu(\mathrm{d}u)\,\mathrm{d}s\} E^Q[L(\tau_n)|\mathcal{F}^O]$ recursively. We recall that when a single quake is observed at time $\tau_i$ it has no impact on the cluster quakes $E_{i+1} = E_i$ and $D_{\tau_i} = D_{\tau_i-}$, while when a cluster quake is observed at time $\tau_i$, there are two possible scenarios: The first case is $D_{\tau_i-} = 0$, so the observation at time $\tau_i$ is a new mother quake $D_{\tau_i} = 1$ and $E_{i+1} = i$. The second case is $D_{\tau_i-} = 1$, so the observation is an offspring quake; hence, $E_{i+1} = E_i$. This observation kills the cluster with probability $p$. In other words, we look at a Bernoulli($p$) random variable $I_i$, which is independent of everything else. If $I_i = 1$, $D_{\tau_i} = 1$; otherwise $D_{\tau_i} = 0$.

For $0 < i < j$,

$$E^Q[L(\tau_j)\mathbf{1}_{\{D_{\tau_j}=0, E_{j+1}=i\}}|\mathcal{F}^O]$$

$$= E^Q[L(\tau_j)\mathbf{1}_{\{D_{\tau_j}=0, E_{j+1}=i\}}\mathbf{1}_{\{\rho_j=0\}}|\mathcal{F}^O] + E^Q[L(\tau_j)\mathbf{1}_{\{D_{\tau_j}=0, E_{j+1}=i\}}\mathbf{1}_{\{\rho_j=1\}}|\mathcal{F}^O]$$

$$= E^Q[L(\tau_j)\mathbf{1}_{\{D_{\tau_{j-1}}=0, E_j=i\}}\mathbf{1}_{\{\rho_j=0\}}|\mathcal{F}^O] + E^Q[L(\tau_j)\mathbf{1}_{\{D_{\tau_{j-1}}=1, E_j=i\}}\mathbf{1}_{\{\rho_j=1\}}\mathbf{1}_{\{I_j=1\}}|\mathcal{F}^O]$$

$$= E^Q\left[L(\tau_{j-1})\gamma_j \exp\left\{-\int_{E \times [\tau_{j-1},\tau_j]} [\varepsilon(u,s) + \gamma(u,s)]\nu(\mathrm{d}u)\,\mathrm{d}s\right\}\right.$$
$$\left. \times \mathbf{1}_{\{D_{\tau_{j-1}}=0, E_j=i\}}\mathbf{1}_{\{\rho_j=0\}}\Big|\mathcal{F}^O\right]$$
$$+ E^Q\left[L(\tau_{j-1})\lambda_{j,i} \exp\left\{-\int_{E \times [\tau_{j-1},\tau_j]} (\lambda(u,s,y_i) + \gamma(u,s))\nu(\mathrm{d}u)\,\mathrm{d}s\right\}\right.$$
$$\left. \times \mathbf{1}_{\{D_{\tau_{j-1}}=1, E_j=i\}}\mathbf{1}_{\{\rho_j=1\}}\mathbf{1}_{\{I_j=1\}}\Big|\mathcal{F}^O\right]$$

$$= \gamma_j \exp\left\{-\int_{E \times [\tau_{j-1},\tau_j]} (\varepsilon(u,s) + \gamma(u,s))\nu(\mathrm{d}u)\,\mathrm{d}s\right\}$$

$$\times E^Q[L(\tau_{j-1})\mathbf{1}_{D_{\{\tau_{j-1}=0,E_j=i\}}}\mathbf{1}_{\{\rho_j=0\}}|\mathcal{F}^O]$$

$$+ \lambda_{j,i}\exp\left\{-\int_{E\times[\tau_{j-1},\tau_j]}(\lambda(u,s,y_i)+\gamma(u,s))\nu(\mathrm{d}u)\,\mathrm{d}s\right\}$$

$$\times E^Q[L(\tau_{j-1})\mathbf{1}_{\{D_{\tau_{j-1}=1,E_j=i}\}}\mathbf{1}_{\{\rho_j=1\}}\mathbf{1}_{\{I_j=1\}}|\mathcal{F}^O]$$

$$=\exp\left\{-\int_{E\times[\tau_{j-1},\tau_j]}\gamma(u,s)\nu(\mathrm{d}u)\,\mathrm{d}s\right\}$$

$$\times\left\{\gamma_j\exp\left[-\int_{E\times[\tau_{j-1},\tau_j]}\varepsilon(u,s)\nu(\mathrm{d}u)\,\mathrm{d}s\right]E^Q[L(\tau_{j-1})\mathbf{1}_{\{D_{\tau_{j-1}=0,E_j=i}\}}\mathbf{1}_{\{\rho_j=0\}}|\mathcal{F}^O]\right.$$

$$+ p\lambda_{j,i}\exp\left[-\int_{E\times[\tau_{j-1},\tau_j]}\lambda(u,s,y_i)\nu(\mathrm{d}u)\,\mathrm{d}s\right]$$

$$\left.\times E^Q[L(\tau_{j-1})\mathbf{1}_{\{D_{\tau_{j-1}=1,E_j=i}\}}\mathbf{1}_{\{\rho_j=1\}}|\mathcal{F}^O]\right\}.$$

The second equality utilizes the fact that the switch $I_j$ is triggered only when a cluster point is observed. The last equality uses the fact that $I_j$ is an independent Bernoulli variable with parameter $p$. Since $\rho_j$ is an independent Bernoulli variable with parameter $1/2$, we have:

$$\exp\left\{\int_{E\times[\tau_{j-1},\tau_j]}\gamma(u,s)\nu(\mathrm{d}u)\,\mathrm{d}s\right\}E^Q[L(\tau_j)\mathbf{1}_{D_{\{\tau_j=0,E_{j+1}=i\}}}|\mathcal{F}^O]$$

$$=\frac{1}{2}\gamma_j\exp\left\{-\int_{E\times[\tau_{j-1},\tau_j]}\varepsilon(u,s)\nu(\mathrm{d}u)\,\mathrm{d}s\right\}E^Q[L(\tau_{j-1})\mathbf{1}_{D_{\{\tau_{j-1}=0,E_j=i\}}}|\mathcal{F}^O]$$

$$+\frac{1}{2}p\lambda_{j,i}\exp\left\{-\int_{E\times[\tau_{j-1},\tau_j]}\lambda(u,s,y_i)\nu(\mathrm{d}u)\,\mathrm{d}s\right\}E^Q[L(\tau_{j-1})\mathbf{1}_{\{D_{\tau_{j-1}=1,E_j=i}\}}|\mathcal{F}^O].$$

We ignore the constant factor $1/2$ in the algorithm.

This raises our forward algorithm. Let $l_j(d,i) \propto E^Q[L(\tau_j)\mathbf{1}_{\{D_{\tau_j=d,E_{j+1}=i}\}}|\mathcal{F}^O]$, where $d=0,1$; $j=1,2,\ldots$; $i=0,\ldots,j-1,j$. ($i=0$ indicates that there is no cluster point.) $l_1(0,0) = \gamma_1\exp\{-\int_{E\times[0,\tau_1]}\varepsilon(u,s)\nu(\mathrm{d}u)\,\mathrm{d}s\}$, $l_1(1,0) = l_1(0,1) = 0$, $l_1(1,1) = \varepsilon_1 \times \exp\{-\int_{E\times[0,\tau_1]}\varepsilon(u,s)\nu(\mathrm{d}u)\,\mathrm{d}s\}$. Furthermore,

$$l_j(0,i) = \begin{cases} \gamma_j\exp\left\{-\int_{E\times[\tau_{j-1},\tau_j]}\varepsilon(u,s)\nu(\mathrm{d}u)\,\mathrm{d}s\right\}l_{j-1}(0,i) \\ \quad + p\lambda_{j,i}\exp\left\{-\int_{E\times[\tau_{j-1},\tau_j]}\lambda(u,s,y_i)\nu(\mathrm{d}u)\,\mathrm{d}s\right\}l_{j-1}(1,i), & 0\leq i<j, \\ 0, & i=j. \end{cases} \quad \text{(A.5)}$$

$$l_j(1,i) = \begin{cases} \gamma_j \exp\left\{-\int_{E\times[\tau_{j-1},\tau_j]} \lambda(u,s,y_i)\nu(\mathrm{d}u)\,\mathrm{d}s\right\} l_{j-1}(1,i) \\ \quad + l_{j-1}(1,i) q \lambda_{j,i} \exp\left\{-\int_{E\times[\tau_{j-1},\tau_j]} \lambda(u,s,y_i)\nu(\mathrm{d}u)\,\mathrm{d}s\right\}, & 0 \le i < j, \\ \sum_{k=0}^{j-1} l_{j-1}(0,k)\varepsilon_j \exp\left\{-\int_{E\times[\tau_{j-1},\tau_j]} \varepsilon(u,s)\nu(\mathrm{d}u)\,\mathrm{d}s\right\}, & i = j. \end{cases} \quad \text{(A.6)}$$

The likelihood until time $\tau_j$ is $L_j = \exp\{-\int_{E\times[0,\tau_j]} \gamma(u,s)\nu(\mathrm{d}u)\,\mathrm{d}s\} \sum_{i=0}^{j} \sum_{d=0}^{1} l_j(d,i)$. Even when there is a moderate number of observations, the scale of the likelihood often exceeds what a computer can handle. Therefore, we compute the $\log(L)$ instead. The trick is normalizing $l_j$ at each step, thus we have the forward algorithm:

1. Computing the normalizing constant: $c_{j-1} = \sum_{i=0}^{j-1} \sum_{d=0}^{1} l_{j-1}(d,i)$.
2. Normalization: $l_{j-1}(d,i) = l_{j-1}(d,i)/c_{j-1}$.
3. Updating $l_j(d,i)$ according to (A.5) and (A.6).
4. $\log(L) = \sum_{j=1}^{n} \log(c_j) - \int_{E\times[0,\tau_n]} \gamma(u,s)\nu(\mathrm{d}u)\,\mathrm{d}s$.

From the discussion above, we can find the likelihood for any specific set of parameters. Looking for the MLE hence is a standard optimization problem. It turns out that a non-derivative method works better in the example of this paper. In particular, we use the Nelder–Mead simplex method to search for the MLE (see [11]).

The asymptotic confidence intervals for the parameters can be constructed. Observe that we have a hidden Markov model (HMM), essentially. Hence the theorems of asymptotic normality of the MLE for a general HMM should apply (see [1, 3, 10]). Consequently, there are corresponding likelihood-ratio tests for the HMM as established in [7]. The asymptotic confidence intervals can then be constructed by inverting the test statistics (see [14]).

### A.3. The most likely cluster sequence

We borrow the Viterbi algorithm from HMM literature to compute the most likely cluster sequence in our setting.

Let $l_j^\star(d,i)$ be the maximum likelihood of all cluster sequences with $D_{\tau_j} = d$ and $y_i$ as the latest mother quake. As in Section A.1,

$$l_1^\star(0,0) = \gamma_1 \exp\left\{-\int_{E\times[0,\tau_1]} \varepsilon(u,s)\nu(\mathrm{d}u)\,\mathrm{d}s\right\}, \qquad l_1^\star(1,0) = l_1^\star(0,1) = 0,$$

$$l_1^\star(1,1) = \varepsilon_1 \exp\left\{-\int_{E\times[0,\tau_1]} \varepsilon(u,s)\nu(\mathrm{d}u)\,\mathrm{d}s\right\},$$

$$l_j^\star(0,i) = \begin{cases} \gamma_j \exp\left\{-\int_{E\times[\tau_{j-1},\tau_j]} \varepsilon(u,s)\nu(\mathrm{d}u)\,\mathrm{d}s\right\} l_{j-1}^\star(0,i), & i = 0, \\ \max\left\{p\lambda_{j,i}\exp\left[-\int_{E\times[\tau_{j-1},\tau_j]} \lambda(u,s,y_i)\nu(\mathrm{d}u)\,\mathrm{d}s\right] l_{j-1}^\star(1,i),\right. \\ \left.\gamma_j \exp\left[-\int_{E\times[\tau_{j-1},\tau_j]} \varepsilon(u,s)\nu(\mathrm{d}u)\,\mathrm{d}s\right] l_{j-1}^\star(0,i)\right\}, & 1 \le i < j, \\ 0, & i = j. \end{cases} \quad (A.7)$$

$$l_j^\star(1,i) = \begin{cases} 0, & i = 0 \\ \max\left\{\gamma_j \exp\left\{-\int_{E\times[\tau_{j-1},\tau_j]} \lambda(u,s,y_i)\nu(\mathrm{d}u)\,\mathrm{d}s\right\} l_{j-1}^\star(1,i),\right. \\ \left.q\lambda_{j,i}\exp\left\{-\int_{E\times[\tau_{j-1},\tau_j]} \lambda(u,s,y_i)\nu(\mathrm{d}u)\,\mathrm{d}s\right\} l_{j-1}^\star(1,i)\right\}, & 1 \le i < j, \\ \max\left\{\varepsilon_j \exp\left\{-\int_{E\times[\tau_{j-1},\tau_j]} \varepsilon(u,s)\nu(\mathrm{d}u)\,\mathrm{d}s\right\} l_{j-1}^\star(0,k);\right. \\ \left.k = 0,1,\ldots,j-1\right\}, & i = j. \end{cases} \quad (A.8)$$

So the procedure to find the most likely cluster sequence starts from the calculation of $l_j^\star(d,i)$, using recursion in (A.7) and (A.8) while always keeping a record of the "winning sequence" in the maximum finding operation. Finally the last state $(d,i)^\star$ is found where

$$(d,i)^\star = \arg\max_{d=0,1;0\le i\le n} l_n^\star(d,i) \quad (A.9)$$

and, starting from this state, the sequence is recovered by backtracking. As before, normalization is necessary in each step of the recursion to prevent them from degenerating to 0 or infinity.

# Appendix B: Proofs

**Proof of Lemma 2.2.**
  By Theorem III.20 of Protter [13]

$$\begin{aligned} M(A,t) &= C(A,t) - \int_{A\times[0,t]} \lambda_Q(u,s)\nu(\mathrm{d}u)\,\mathrm{d}s \\ &\quad - \int_{A\times[0,t]} \frac{1}{L(s)}\left(\frac{\lambda(u,s,D_{s-},\eta_{s-})}{\lambda_Q(u,s)} - 1\right) L(s-)C(\mathrm{d}u\times\mathrm{d}s) \\ &= C(A,t) - \int_{A\times[0,t]} \lambda_Q(u,s)\nu(\mathrm{d}u)\,\mathrm{d}s \\ &\quad - \int_{A\times[0,t]} \frac{\lambda(u,s,D_{s-},\eta_{s-}) - \lambda_Q(u,s)}{\lambda(u,s,D_{s-},\eta_{s-})} C(\mathrm{d}u\times\mathrm{d}s) \end{aligned}$$

is a local martingale, and hence

$$\int_{A\times[0,t]} \frac{\lambda(u,s,D_{s-},\eta_{s-})}{\lambda_Q(u,s)} M(\mathrm{d}u \times \mathrm{d}s) = C(A,t) - \int_{A\times[0,t]} \lambda(u,s,D_{s-},\eta_{s-})\nu(\mathrm{d}u)\,\mathrm{d}s$$

is as well. $\square$

**Proof of Theorem 2.4.** To simplify the notation, we use $f(0,\eta_{s-}+\delta_{(u,s)})$ to denote $f(\cdot+\delta_{(0,s)},\eta_{s-}+\delta_{(u,s)})$ and $f(1,\eta_{s-}+\delta_{(u,s)})$ to denote $f(\cdot+\delta_{(1,s)},\eta_{s-}+\delta_{(u,s)})$.

Noting that

$$[f,L]_t = \int_{E\times[0,t]} \{\mathbf{1}_{\{D_{s-}=0\}}[f(1,\eta_{s-}+\delta_{(u,s)}) - f(h_{s-},\eta_{s-})]$$
$$+ \mathbf{1}_{\{D_{s-}=1\}}[f((1-I_{C(E,s)})D_{s-},\eta_{s-}+\delta_{(u,s)}) - f(h_{s-},\eta_{s-})]\}$$
$$\times \left(\frac{\lambda(u,s,\eta_{s-})}{\lambda_Q(u,s)} - 1\right) L(s-)C(\mathrm{d}u \times \mathrm{d}s),$$

$$f(h_t,\eta_t)L(t) = f(h_0,\eta_0) + \int_0^t f(h_{s-},\eta_{s-})\,\mathrm{d}L(s)$$
$$+ \int_0^t L(s-)\,\mathrm{d}f\circ(h(s),\eta(s)) + [f,L]_t$$
$$= f(h_0,\eta_0) + \int_0^t f(h_{s-},\eta_{s-})\,\mathrm{d}L(s)$$
$$+ \int_{E\times[0,t]} \{\mathbf{1}_{\{D_{s-}=0\}}[f(1,\eta_{s-}+\delta_{(u,s)}) - f(h_{s-},\eta_{s-})]$$
$$+ \mathbf{1}_{\{D_{s-}=1\}}[f(1-I_{C(E,s)},\eta_{s-}+\delta_{(u,s)}) - f(h_{s-},\eta_{s-})]\}$$
$$\times \frac{\lambda(u,s,\eta_{s-})}{\lambda_Q(u,s)}L(s-)C(\mathrm{d}u \times \mathrm{d}s)$$
$$= f(h_0,\eta_0) + \int_{E\times[0,t]} f(h_{s-},\eta_{s-})\left[\frac{\lambda(u,s,h_{s-},\eta_{s-})}{\lambda_Q(u,s)} - 1\right]L(s-)$$
$$\times [C(\mathrm{d}u \times \mathrm{d}s) - \lambda_Q(u,s)\nu(\mathrm{d}u)\,\mathrm{d}s]$$
$$+ \int_{E\times[0,t]} \{\mathbf{1}_{\{D_{s-}=0\}}[f(1,\eta_{s-}+\delta_{(u,s)}) - f(h_{s-},\eta_{s-})]$$
$$+ \mathbf{1}_{\{D_{s-}=1\}}[f((1-I_{C(E,s)})D_{s-},\eta_{s-}+\delta_{(u,s)}) - f(h_{s-},\eta_{s-})]\}$$
$$\times \frac{\lambda(u,s,h_{s-},\eta_{s-})}{\lambda_Q(u,s)}L(s-)C(\mathrm{d}u \times \mathrm{d}s)$$

$$= f(h_0, \eta_0) - \int_{E \times [0,t]} f(h_s, \eta_s)[\lambda(u, s, h_s, \eta_s) - \lambda_Q(u, s)]L(s-)\nu(\mathrm{d}u)\,\mathrm{d}s$$

$$+ \int_{E \times [0,t]} \{[\mathbf{1}_{\{D_{s-}=0\}} f(1, \eta_{s-} + \delta_{(u,s)})$$

$$+ \mathbf{1}_{\{D_{s-}=1\}} f(1 - I_{C(E,s)}, \eta_{s-} + \delta_{(u,s)})] \frac{\lambda(u, s, h_{s-}, \eta_{s-})}{\lambda_Q(u, s)}$$

$$- f(h_{s-}, \eta_{s-})\}L(s-)\rho_{O(E,s)}(u)O(\mathrm{d}u \times \mathrm{d}s),$$

where

$$C(A, t) = \int_{A \times [0,t]} \rho_{O(E,s)}(u)O(\mathrm{d}u \times \mathrm{d}s)$$

and under the reference measure, the $\{\rho_k(\cdot), k = 1, 2, \ldots\}$ are independent with $Q\{\rho_k(u) = 1\} = 1 - Q\{\rho_k(u) = 0\} = \lambda_Q(u, s)/[\lambda_Q(u, s) + \gamma(u, s)]$ and are independent of $O$.

Averaging out the random variables that are independent of $O$ under $Q$, the equation for the unnormalized conditional distribution becomes

$$\phi(f, t) = \phi(f, 0) - \int_{E \times [0,t]} \phi(f(\cdot, \cdot)[\lambda(u, s, \cdot, \cdot) - \lambda_Q(u, s)], s)\nu(\mathrm{d}u)\,\mathrm{d}s$$

$$+ \int_{E \times [0,t]} \phi\left(f_{\mathrm{new}} \frac{\lambda(u, s, \cdot, \cdot)}{\lambda_Q(u, s)} - f(\cdot, \cdot), s-\right) \frac{\lambda_Q(u, s)}{\lambda_Q(u, s) + \gamma(u, s)} O(\mathrm{d}u \times \mathrm{d}s).$$

Applying Itô's formula,

$$\pi(f, t) = \pi(f, 0)$$

$$+ \int_{E \times [0,t]} \frac{\pi(f_{\mathrm{new}} \lambda(u, s, \cdot, \cdot), s-) - \pi(\lambda(u, s, \cdot, \cdot), s-)\pi(f, s-)}{\pi(\lambda(u, s, \cdot, \cdot) + \gamma(u, s), s-)} O(\mathrm{d}u \times \mathrm{d}s)$$

$$- \int_{E \times [0,t]} \pi(f(\cdot, \cdot)\lambda(u, s, \cdot, \cdot), s) - \pi(f, s)\pi(\lambda(u, s, \cdot, \cdot), s)\nu(\mathrm{d}u)\,\mathrm{d}s. \qquad \square$$

**Proof of Theorem 3.3.** We apply Theorem 2.4. For $\tau_k \leq t < \tau_{k+1}$,

$$\pi(\theta_0(y)\alpha, t) = \pi(\theta_0(y)\alpha, \tau_k) - \int_{E \times [\tau_k, t]} \pi(\theta_0(y)\alpha\lambda(u, s, h, \eta), s)\nu(\mathrm{d}u)\,\mathrm{d}s$$

$$+ \int_{E \times [\tau_k, t]} \pi(\theta_0(y)\alpha, s)\pi(\lambda(u, s, h.\eta), s)\nu(\mathrm{d}u)\,\mathrm{d}s$$

$$= \pi(\theta_0(y)\alpha, \tau_k) - \int_{E \times [\tau_k, t]} \pi\left(\sum_{x \in O(s)} \theta_0(y)\theta_0(x)\alpha\lambda(u, s, x), s\right)\nu(\mathrm{d}u)\,\mathrm{d}s$$

$$+ \int_{E \times [\tau_k, t]} \pi(\theta_0(y)\alpha, s)$$

$$\times \left[\sum_{x \in O(s)} \pi(\theta_0(x)\alpha(\lambda(u,s,x) - \varepsilon(u,s)), s) + \varepsilon(u,s)\right] \nu(\mathrm{d}u)\,\mathrm{d}s$$

$$= \pi(\theta_0(y)\alpha, \tau_k) - \int_{[\tau_k, t]} \sum_{x \in O(s)} \pi(\theta_0(y)\theta_0(x)\alpha, s) \int_E \lambda(u,s,x)\nu(\mathrm{d}u)\,\mathrm{d}s$$

$$+ \int_{[\tau_k, t]} \pi(\theta_0(y)\alpha, s)$$

$$\times \left[\sum_{x \in O(s)} \pi(\theta_0(x)\alpha, s) \int_E \lambda(u,s,x) - \varepsilon(u,s)\nu(\mathrm{d}u) + \varepsilon(s)\right] \mathrm{d}s$$

$$= \pi(\theta_0(y)\alpha, \tau_k) - \int_{[\tau_k, t]} \pi(\theta_0(y)\alpha, s) a(y, s) \mathbf{1}_{\{y \in O(\tau_k)\}}\,\mathrm{d}s$$

$$+ \int_{[\tau_k, t]} \pi(\theta_0(y)\alpha, s) \left[\sum_{x \in Y(\tau_k)} \pi(\theta_0(x)\alpha, s)(a(x,s) - \varepsilon(s)) + \varepsilon(s)\right] \mathrm{d}s.$$

We use the fact that $\alpha = \sum_{x \in O(s)} \theta_0(x)\alpha$ to get the second equality. Thus, for $i = 1, 2, \ldots, k$,

$$\frac{\mathrm{d}\pi(\theta_0(y_i), t)}{\mathrm{d}t} = \left\{-a(y_i, t) + \varepsilon(s) + \sum_{j=1}^{k} [a(y_j, t) - \varepsilon(s)]\pi(\theta_0(y_j), t)\right\} \pi(\theta_0(y_i), t), \quad \text{(B.1)}$$

(3.2) is the unique solution of this system of ordinary differential equations.

At time $\tau_{k+1}$,

$$\pi(\theta_0(y)\alpha, \tau_{k+1}) = \pi(\theta_0(y)\alpha, \tau_{k+1}-)$$

$$+ \frac{\pi([1-\alpha(\cdot,\cdot)][\theta_0(y)(1,\cdot+\delta_{(u,s)})]\lambda(y_{k+1}, \tau_{k+1}, \cdot, \cdot), \tau_{k+1}-)}{d_{k+1}}$$

$$+ \frac{\pi(\alpha(\cdot,\cdot)q[(\theta_0(y)\alpha)(1,\cdot+\delta_{(u,s)})]\lambda(y_{k+1}, \tau_{k+1}, \cdot, \cdot), \tau_{k+1}-)}{d_{k+1}}$$

$$- \frac{\pi(\lambda(y_{k+1}, \tau_{k+1}, \cdot, \cdot), \tau_{k+1}-)\pi(\theta_0(y)\alpha, \tau_{k+1}-)}{d_{k+1}}$$

$$= \frac{\pi(\theta_0(y)\alpha, \tau_{k+1}-)\gamma_{k+1}}{d_{k+1}}$$

$$+ \frac{\pi([1-\alpha(\cdot,\cdot)][\theta_0(y)(1,\cdot+\delta_{(u,s)})]\lambda(y_{k+1}, \tau_{k+1}, \cdot, \cdot), \tau_{k+1}-)}{d_{k+1}}$$

$$+ \frac{\pi(\alpha(\cdot,\cdot)q[(\theta_0(y)\alpha)(1,\cdot+\delta_{(u,s)})]\lambda(y_{k+1}, \tau_{k+1}, \cdot, \cdot), \tau_{k+1}-)}{d_{k+1}}.$$

For $i < k+1$,

$$\pi(\theta_0(y_i)\alpha, \tau_{k+1}) = \frac{\pi(\theta_0(y)\alpha, \tau_{k+1}-)(\gamma_{k+1} + q\lambda_{k+1,i})}{d_{k+1}}.$$

For $i = k+1$,

$$\pi(\theta_0(y_{k+1})\alpha, \tau_{k+1}) = \frac{\sum_{j=1}^{k}(q\lambda_{k+1,j} - \varepsilon_{k+1})\pi(\theta_0(y_j)\alpha, \tau_{k+1}-) + \varepsilon_{k+1}}{d_{k+1}}. \qquad \square$$

**Proof of Theorem 3.2.** For $\tau_k \leq t < \tau_{k+1}$,

$$\pi(\theta_0(y)\alpha f, t) = \pi(\theta_0(y)\alpha f, \tau_k) - \int_{E\times[\tau_k,t]} \pi(\theta_0(y)\alpha f\lambda(u,s,h,\eta), s)\nu(\mathrm{d}u)\,\mathrm{d}s$$

$$+ \int_{E\times[\tau_k,t]} \pi(\theta_0(y)\alpha f, s)\pi(\lambda(u,s,h,\eta), s)\nu(\mathrm{d}u)\,\mathrm{d}s$$

$$= \pi(\theta_0(y)\alpha f, \tau_k) - \int_{E\times[\tau_k,t]} \pi\left(\sum_{x\in O(s)} \theta_0(y)\theta_0(x)\alpha f\lambda(u,s,x), s\right)\nu(\mathrm{d}u)\,\mathrm{d}s$$

$$+ \int_{E\times[\tau_k,t]} \pi(\theta_0(y)\alpha f, s)$$

$$\times \left\{\sum_{x\in O(s)} \pi(\theta_0(x)\alpha[\lambda(u,s,x) - \varepsilon(u,s)], s) + \varepsilon(u,s)\right\}\nu(\mathrm{d}u)\,\mathrm{d}s$$

$$= \pi(\theta_0(y)\alpha f, \tau_k) - \int_{[\tau_k,t]} \sum_{x\in O(s)} \pi(\theta_0(y)\theta_0(x)\alpha f, s)\int_E \lambda(u,s,x)\nu(\mathrm{d}u)\,\mathrm{d}s$$

$$+ \int_{[\tau_k,t]} \pi(\theta_0(y)\alpha f, s)$$

$$\times \left\{\sum_{x\in O(s)} \pi(\theta_0(x)\alpha, s)\int_E \lambda(u,s,x) - \varepsilon(u,s)\nu(\mathrm{d}u) + \varepsilon(s)\right\}\mathrm{d}s$$

$$= \pi(\theta_0(y)\alpha f, \tau_k) - \int_{[\tau_k,t]} \pi(\theta_0(y)\alpha f, s)a(y,s)\mathbf{1}_{\{y\in O(\tau_k)\}}\,\mathrm{d}s$$

$$+ \int_{[\tau_k,t]} \pi(\theta_0(y)\alpha f, s)$$

$$\times \left\{\sum_{x\in Y(\tau_k)} \pi(\theta_0(x)\alpha, s)[a(x,s) - \varepsilon(s)] + \varepsilon(s)\right\}\mathrm{d}s.$$

Thus,

$$\frac{\mathrm{d}\pi(\theta_0(y_i)\alpha f, t)}{\mathrm{d}t}$$
$$= \left\{-a(y_i,t) + \varepsilon(s) + \sum_{j=1}^{k}[a(y_j,t) - \varepsilon(s)]\pi(\theta_0(y_j),t)\right\}\pi(\theta_0(y_i)\alpha f, t), \quad \text{(B.2)}$$
$$i = 1, 2, \ldots, k$$

Comparing (B.2) and (B.1), we conclude:

$$\pi(\theta_0(y)\alpha f, s) = \frac{\pi(\theta_0(y)\alpha, s)\pi(\theta_0(y)\alpha f, \tau_k)}{\pi(\theta_0(y)\alpha, \tau_k)}$$

and (3.3) follows easily.

$$\pi(\theta_0(y)\alpha f, \tau_{k+1}) = \pi(\theta_0(y)\alpha f, \tau_{k+1}-) - \frac{\pi(\lambda(y_{k+1}, \tau_{k+1}, \cdot, \cdot), \tau_{k+1}-)\pi(\theta_0(y)\alpha f, \tau_{k+1}-)}{d_{k+1}}$$
$$+ \frac{\pi([1-\alpha(\cdot,\cdot)][(\theta_0(y)f)(1,\cdot+\delta_{(u,s)})]\lambda(y_{k+1},\tau_{k+1},\cdot,\cdot), \tau_{k+1}-)}{d_{k+1}}$$
$$+ \frac{\pi(\alpha(\cdot,\cdot)q[(\theta_0(y)f)(1,\cdot+\delta_{(u,s)})]\lambda(y_{k+1},\tau_{k+1},\cdot,\cdot), \tau_{k+1}-)}{d_{k+1}}$$
$$= \frac{[1 - \delta_{y_{k+1}}(y)]\pi(\theta_0(y)\alpha f, \tau_{k+1}-)\gamma_{k+1}}{d_{k+1}}$$
$$+ \frac{[1 - \delta_{y_{k+1}}(y)]q\sum_{j=1}^{k}\lambda_{k+1,j}\pi(f(1,\cdot+\delta_{y_{k+1}})\theta_0(y_j)\alpha, \tau_{k+1}-)}{d_{k+1}}$$
$$+ \frac{\delta_{y_{k+1}}(y)\pi(f(1,\cdot+\delta_{y_{k+1}})(1-\alpha), \tau_{k+1}-)\varepsilon_{k+1}}{d_{k+1}},$$

$$\pi(f, \tau_{k+1}-)$$
$$= \pi(f, \tau_k) - \int_{E\times[\tau_k,\tau_{k+1}]} \pi(f\lambda(u,s,h,\eta), s)\nu(\mathrm{d}u)\,\mathrm{d}s$$
$$+ \int_{E\times[\tau_k,\tau_{k+1}]} \pi(f, s)\pi(\lambda(u,s,h,\eta), s)\nu(\mathrm{d}u)\,\mathrm{d}s$$
$$= \pi(f, \tau_k) - \int_{E\times[\tau_k,\tau_{k+1}]} \pi\left(\sum_{x\in O(s)}\theta_0(x)\alpha f\lambda(u,s,x) + (1-\alpha)\varepsilon(u,s)f, s\right)\nu(\mathrm{d}u)\,\mathrm{d}s$$
$$+ \int_{E\times[\tau_k,\tau_{k+1}]} \pi(f, s)\left[\sum_{x\in O(s)}\pi(\theta_0(x)\alpha\lambda(u,s,x), s) + (1-\alpha)\varepsilon(u,s)\right]\nu(\mathrm{d}u)\,\mathrm{d}s$$

$$= \pi(f, \tau_k) - \int_{[\tau_k, \tau_{k+1}]} \sum_{x \in O(s)} \pi(\theta_0(x)\alpha f, s)[a(x,s) - \varepsilon(s)] \, \mathrm{d}s$$

$$+ \int_{[\tau_k, \tau_{k+1}]} \pi(f, s) \sum_{x \in O(\tau_k)} \pi(\theta_0(x)\alpha, s)[a(x,s) - \varepsilon(s)] \, \mathrm{d}s,$$

where $a(x,s) = \int_E \lambda(u, s, x)\nu(\mathrm{d}u)$, $\varepsilon(s) = \int_E \varepsilon(u,s)\nu(\mathrm{d}u)$. Hence,

$$\frac{\mathrm{d}\pi(f,t)}{\mathrm{d}t} = -\sum_{x \in Y(t)} \pi(\theta_0(x)\alpha f, t)[a(x,t) - \varepsilon(s)]$$

$$+ \pi(f,t) \left\{ \sum_{x \in Y(\tau_k)} \pi(\theta_0(x)\alpha, t)[a(x,t) - \varepsilon(s)] \right\},$$

$$\pi(f, \tau_{k+1}) = \pi(f, \tau_{k+1}-) - \frac{\pi(\lambda(y_{k+1}, \tau_{k+1}, \cdot, \cdot), \tau_{k+1}-)\pi(\alpha f, \tau_{k+1}-)}{d_k + 1}$$

$$+ \frac{\pi(f_{\mathrm{new}}\lambda(y_{k+1}, \tau_{k+1}, \cdot, \cdot), \tau_{k+1}-)}{d_{k+1}}$$

$$= \frac{\pi(f, \tau_{k+1}-)\gamma_{k+1} + \pi(f(1, \cdot + \delta_{y_{k+1}})(1-\alpha)\varepsilon_{k+1}, \tau_{k+1}-)}{d_{k+1}}$$

$$+ \frac{\sum_{j=1}^{k} \lambda_{k+1,j} \pi([pf(0, \cdot + \delta_{y_{k+1}}) + qf(1, \cdot + \delta_{y_{k+1}})]\theta_0(y_j)\alpha, \tau_{k+1}-)}{d_{k+1}}. \qquad \square$$

**Proof of Theorem A.1.** Apply Lemma 2.2 and note the independence of $C$ and $N$. $\square$

## Acknowledgements

This material is based upon work supported by, or in part by, the U.S. Army Research Laboratory and the U.S. Army Research Office under contract, grant number DAAD19-01-1-0502, and by NSF Grant DMS 05-03983.

I am deeply grateful to my adviser Tom Kurtz for very helpful discussions and comments. He taught me the filtering idea and directed my work on its applications. We thank Jiancang Zhuang who kindly sent us the data set and his numerical results, which were plotted in [8]. We also thank Feng-Chang Lin for pointing us to Zhuang's paper.